# Anisotropic connections and parallel transport in Finsler spacetimes[★]

Miguel Ángel Javaloyes[1], Miguel Sánchez[2] and Fidel F. Villaseñor[2]

[1] University of Murcia, Department of Mathematics, Spain
`majava@um.es`,
[2] Universidad de Granada, Departmento de Geometría y Topología, Spain
`sanchezm@ugr.es, fidelfv@ugr.es`

**Abstract.** The general notion of anisotropic connections $\nabla$ is revisited, including its precise relations with the standard setting of pseudo-Finsler metrics, i.e., the canonic nonlinear connection and the (linear) Finslerian connections. In particular, the vertically trivial Finsler connections are identified canonically with anisotropic connections. So, these connections provide a simple intrinsic interpretation of a part of any Finsler connection closer to the Koszul formulation in $M$. Moreover, a new covariant derivative and parallel transport along curves is introduced, taking first a self-propagated vector (*instantaneous observer*) so that it serves as a reference for the propagation of the others. The covariant derivative of any anisotropic tensor is given by the natural derivative of a curve of tensors obtained by parallel transport along a curve and, in the case of pseudo-Finsler metrics, this is used to characterize the Levi-Civita–Chern anisotropic connection as the one that preserves the length of parallely propagated vectors.

## 1 Introduction

The standard geometric picture for the description of a (pseudo) Finsler metric $L = F^2 : A(\subset TM) \to \mathbb{R}$ comprises two elements: a nonlinear connection $\nu \equiv N_i^a(x, y)$ on the fibration $A \to M$ and a linear connection $\nabla^* \equiv (H_{ij}^k, V_{ij}^k)$ on the vertical vector bundle $\mathcal{V}A \to A$. The former is canonically associated with the geodesics of $L$, which are regarded as critical points of the energy functional. However, there are quite a few non-equivalent choices for the latter (Berwald, Cartan, Chern, Hashiguchi...). Motivated by the complexity of this and other settings, some researchers have introduced the concept of *anisotropic connection*, a natural generalization of the (pseudo) Riemannian setting which incorporates in a natural way the geometric direction-dependence of Finslerian geometries, [23,24].

Recently, one of the authors has developed systematically the *anisotropic calculus* [10,11], namely, how to make computations with an anisotropic connection, which can be seen as a natural and intuitive generalization of the usual Koszul connections. Some applications have been obtained in [9,16]. In the present article, we revisit this notion, showing precisely its relations with the other elements of the standard setting and providing a further insight on its associated parallel transport.

More precisely, in §2 we introduce heuristically the notions of pseudo-Finsler metric, by looking for general ways of measuring the lengths of curves, and Finsler spacetime, by stressing geometric elements related with measurements. In §3 anisotropic tensor fields on $M$ are introduced and the concept of anisotropic connection $\nabla$ is defined. First, $\nabla$ is regarded as a type of derivative which applies to usual vector fields $X, Y$ on $M$ so that it provides an anisotropic vector field $\nabla_X Y$. Then, the usual rules of derivations are discussed so that $\nabla_X$ can be applied to any anisotropic tensor. We emphasize some issues which will become relevant later such as homogeneity (natural invariance by homotheties), the notion of torsion or the affine structure of the space of all the anisotropic connections.

In §4 and §5 we make a detailed study of the relation between anisotropic connections $\nabla$ and, resp., nonlinear connections $\nu$ on $A \to M$ and linear connections $\nabla^*$ on $\mathcal{V}A \to A$ ($\nu$ and $\nabla^*$ are not assumed to come from any Finsler $L$ a priori). A detailed correspondence is given for the former, in particular:

> Any anisotropic connection $\nabla$ with Christoffel symbols $\Gamma_{ij}^a$ is characterized by a pair composed by a nonlinear connection $N_i^a = \Gamma_{ij}^a y^j$ and a tensor $Q$ satisfiying $Q_{ij}^a y^j = 0$. In the homogeneous case, all the nonlinear connections can be obtained from anisotropic ones (Theorem 2).

The relation between anisotropic $\nabla$ and linear $\nabla^*$ becomes subtler. Indeed, if we are given an auxiliary nonlinear connection $\overset{o}{\nu}$, then $\nabla^*$ can be determined by specifying the covariant derivatives (of the sections of $\mathcal{V}A \to A$) with respect to the horizontal and vertical directions of $\overset{o}{\nu}$. This is standard in Finsler Geometry, and the connections with vanishing vertical derivatives are called *vertically trivial* here; clearly, they are independent of $\overset{o}{\nu}$. Such trivial connections can be put in one to one correspondence with the anisotropic connections by using $\overset{o}{\nu}$. However, this correspondence also becomes independent of $\overset{o}{\nu}$, summing up:

> There is a natural bijection between vertically trivial connections $\nabla^*$ and anisotropic connections $\nabla$, which maps (also bijectively) invariant by homotheties vertically trivial $\nabla^*$'s into homogeneous $\nabla$'s (Proposition 2, Theorem 3).

In §6 we focus on the pseudo-Finsler case. As explained therein, the last result above becomes essential for the identification of anisotropic connections in the pseudo-Finsler setting. Indeed, the 2-homogeneity of $L$ leads to the homogeneity of the involved linear connection $\nabla^*$ (and the canonic nonlinear one $\overset{o}{\nu}$). Some Finsler connections such as Berwald or Chern are vertically trivial and, thus, directly identifiable with anisotropic connections. Moreover, the non vertically trivial ones, as Cartan or Hasiguchi, will project on vertically trivial ones (by using $\overset{o}{\nu}$). So, anisotropic connections provide the non-vertical part of any Finslerian connection, expressed tidily as Koszul-type derivations on $M$. As already pointed out in [10], the metric $L$ allows one to select a unique Levi-Civita anisotropic connection, which is then identifiable to Chern's.

In §7, we introduce the covariant derivative $D_\gamma$ and parallel transport along curves $\gamma$ for any anisotropic connection $\nabla$. Taking into account the dependence of $\nabla$ on the direction, one can choose a reference $W$ (a vector field on $\gamma$ which takes values on $A \subset M$) as in [2, page 121], to define its associated covariant derivative $D_\gamma^W$ and $W$-parallel transport, which will behave as the usual (isotropic) one. This parallel transport is of crucial importance, as it can be used to define in a very natural way the covariant derivative of tensors. As in Finsler Geometry there is a dependence on the direction, one has to be cautious. As a first step, one can parallel transport the observer, which in Finsler spacetimes is interpreted as the direction on the tangent bundle where we are doing the computations. This is defined using a parallel observer determined by $D_\gamma^V V = 0$, a nonlinear equation whose solutions may not be extended on the whole $\gamma$. However, they do extend in the most interesting cases, such as the standard Finsler and the Lorentz-Finsler ones. Once we have a parallel observer along the curve, we can make the parallel transport of any other vector using as a reference this parallel observer. The parallel transport of the observer coincides with the one provided by a non-linear connection (see for example [17, Ch. VII],[1, §2.1.6], [24, page 103], [4, §2.1], [7, Def. 1.4] and [26, §7.6]). However, as far as we know, the second parallel transport with respect to an observer has not been considered in literature.

It turns out that the most economical way to codify all the information of the covariant derivatives along curves in a smooth setting with natural assumptions is with an *anisotropic connection*, which allow for covariant derivatives of any kind of tensor (see Theorem 1). These covariant derivatives were introduced in [11] from a rather abstract viewpoint as tensor derivations which satisfy the Leibniz rule of the tensor product and commute with contractions. To enhance the geometric meaning of these covariant derivatives, we will show in Theorem 5 that they coincide with the (usual) derivative in a vector space of the curve of tensors obtained with parallel transport (and therefore using only covariant derivatives along curves). Finally, one can wonder if given a pseudo-Finsler metric, there is an anisotropic connection with a parallel transport which preserves the length of vectors. We show in §7.3 that this connection exists and it is the Levi-Civita-Chern anisotropic connection, which can be identified (as a vertically trivial linear connection) with the classical Chern connection.

Finally, we would like to emphasize that our approach is useful in the classical Finsler case a well as its (positive definite) variants, such as Kropina, Randers-Kropina, conic and wind Finsler metrics [5,8,13,15,28]. In this article, we emphasize the case of Finsler spacetimes (introducing notions such as *observer*) not only because these constitute an active topic of research where our setting applies naturally [3,8,19], but also because its physical intuitions suggest interesting geometric definitions valid even for the positive definite case.

## 2     General background

### 2.1     Pseudo-Finsler metrics

Let us pose the following problem. Given a manifold $M$, we want to define a general smooth structure that allows us to measure the length of curves. It seems quite natural that this length should be defined as

$$\ell(\gamma) := \int_a^b F(\dot{\gamma})ds$$

for some function $F : TM \to \mathbb{R}$ in the tangent bundle $TM$. From a geometrical viewpoint, this definition should not depend on the parametrization of $\gamma$, which can be achieved by requiring that $F$ is a homogeneous function of degree 1 when restricted to every tangent space. Moreover, if one wants to include relativistic measures and to remain in the smooth realm, it is better to consider the square $L = F^2$, because otherwise one would find many examples where $F$ is not smooth on lightlike vectors, namely, vectors $v \in TM$ where $F(v) = 0$. Indeed, this happens when one considers the one-homogeneous function $F : \mathbb{R}^4 \to \mathbb{R}$ given by

$$F(\tau, v^1, v^2, v^3) = \sqrt{\tau^2 - (v^1)^2 - (v^2)^2 - (v^3)^2},$$

which is non-smooth on the lightlike vectors.

We will make two additional assumptions on $L$.

(i) The first one is that $L$ is not necessarily defined in the whole tangent bundle $TM$, but only in some directions. Sometimes, there are some forbidden directions because of some constraints of the problem, or as in General Relativity, because only trajectories with directions on a cone (say, the future-directed timelike one) will become relevant. Therefore, we will choose as domain of $L$ a subset $A$ in $TM$ which is conic, to permit arbitrary positive reparametrizations of the curves, and open for the sake of simplicity, even though the boundary of $A$ can be considered in different ways (as in wind Riemannian metrics [5] or Finsler spacetimes [14]).

(ii) The second one is related to the (vertical) Hessian of $L$,

$$g_v(u, w) := \frac{1}{2} \frac{\partial^2}{\partial t \partial s} L(v + tu + sw)|_{t=s=0},$$

which, as we will see later, can be thought as the best scalar product approximation of $L$. We will assume that this scalar product is nondegenerate for every $v \in A$ but not necessarily Euclidean (positive definite). Nondegenericity will be essential to obtain the existence and uniqueness of the covariant derivative.

Summing up, the following notion of *pseudo-Finsler metric* collects all the conditions above for a very general definition of length of curves.

**Definition 1.** *Let $M$ be an $n$-manifold, $\pi : TM \to M$ the natural projection of $TM$ onto $M$ and $A \subset TM \setminus \mathbf{0}$ an open subset of $TM$ which is conic, (namely, for every $v \in A$ and $\lambda > 0$, $\lambda v \in A$) and satisfies $\pi(A) = M$. A smooth function $L : A \to \mathbb{R}$ is a* pseudo-Finsler metric *if*

*(1) $L$ is positive homogeneous of degree 2, that is, $L(\lambda v) = \lambda^2 L(v)$ for every $v \in A$ and $\lambda > 0$,*

*(2) for every $v \in A$, the fundamental tensor of $L$, defined as*

$$g_v(u, w) := \frac{1}{2} \frac{\partial^2}{\partial t \partial s} L(v + tu + sw)|_{t=s=0}$$

*for any $u, w \in T_{\pi(v)}M$, is nondegenerate.*

In this definition we have excluded the zero section from $A$. As $A$ is open and conic, the only case in which the zero section could be contained in $A$ is when it is the whole tangent bundle. But even in this case, there are problems with the zero section, because $L$ can be $C^2$ on the zero section only if it comes pointwise from a scalar product (indeed, if $g$ is one half the Hessian of $L$ at the 0 vector of each tangent space, then $L(v) = g(v, v)$ for every $v \in TM$, see [27, Proposition 4.1]).

Given a pseudo-Finsler metric $L : A \to \mathbb{R}$ on a manifold $M$, for every $p \in M$, we define the *indicatrix* at $p$ as

$$\Sigma_p = \{v \in T_pM \cap A : L(v) = 1\},$$

(sometimes the indicatrix of $-L$ may be of interest too) and the *lightcone* as

$$\mathcal{C}_p = \{v \in T_pM \cap A : L(v) = 0\}.$$

Given $p \in M$ and $v \in T_pM$, let us discuss why $g_v$ is the best scalar product approximation of $L$ at $v \in A$. Assume for example that $L(v) = 1$, which can be assumed by homogeneity if $L(v) > 0$ (the other case is analogous). Recall that the restriction

$$g_v|_{\Sigma} : T_v\Sigma_p \times T_v\Sigma_p \to \mathbb{R}$$

coincides with second fundamental form of $\Sigma_p$ with respect to the opposite of the position vector $v$ computed with the affine connection of $T_pM$ (see for example [13, Eq. (2.3)]). Moreover, one has that $v$ is $g_v$-orthogonal to $T_v\Sigma_p$ and, by homogeneity (applying Euler's Theorem), that $g_v(v, v) = L(v)$. This implies that

$$\Sigma^{g_v} = \{w \in T_pM : g_v(w, w) = 1\}$$

satisfies that $T_v\Sigma^{g_v} = T_v\Sigma_p$ and the second fundamental form of $\Sigma^{g_v}$ at $v$ with respect to the opposite of the position vector $v$ coincides with that of $\Sigma_p$.

## 2.2   Finsler spacetimes and its restspace

To generalize the definition of spacetime in a certain manifold $M$, the following observations are in order:

(i)   We need to measure the length of curves to obtain the elapsed time along the trajectory. By the discussion in the previous section, this leads us to consider a pseudo-Finsler metric $L : A \subset TM \setminus 0 \to \mathbb{R}$.

(ii)   Locally, it must approximate the Lorentz-Minkowski spacetime. This implies that for every $v \in A$, the scalar product $g_v$ must be of Lorentz type since, as argued above, $g_v$ is the best approximation of $L$ around $v$.

(iii)   There have to be some vectors with zero length, which are the directions of light rays.

(iv)   Moreover, these lightlike directions must be the limit of the timelike directions, therefore, their boundary.

**Definition 2.** *A* Finsler spacetime *is an $n$-manifold $M$, $n \geq 2$, endowed with a pseudo-Finsler metric $L : A \to (0, +\infty)$ such that*

(i)   *$L$ is a Lorentz-Finsler metric (its indicatrix is strongly concave or equivalently the index of $g_v$ is $n-1$).*

(ii)   *$L$ extends as zero to the closure $\bar{A}$ of $A$ in $TM \setminus 0$ and this extension is smooth with nondegenerate $g_v$.*

(iii)   *For every $p \in M$, $A_p := A \cap T_pM$ is connected, convex and salient (in fact, the last two follow from the other hypotheses, see [14, Remark 3.6]).*

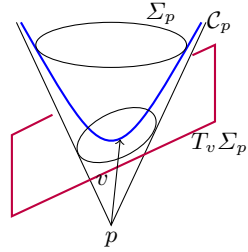

Moreover, the future-directed timelike unit vectors of the indicatrix $\Sigma_p = \{v \in T_pM : L(v) = 1\}$ are used to model the instantaneous observers, while the vectors in the null cone $\mathcal{C}_p = \{v \in T_pM : L(v) = 0\}$ are the lightlike future-directed vectors. The tangent space $T_v\Sigma_p = \{w \in T_pM : g_v(v,w) = 0\}$ is interpreted as the instantaneous restspace of $v$. Even if we assume that $L$ is defined only in $\bar{A}$, it is possible to extend $L$ (in a non-unique way) to the whole tangent bundle (see [19]). For more details about the interpretation of the restspace see [3].

An *observer* in a Finsler spacetime is a unit (future-directed) timelike curve $\gamma : I = (a,b) \to (M,L)$, namely, $\dot{\gamma}(s) \in A$ and $L(\dot{\gamma}(s)) = 1$ for all $s \in I$. The *(instantaneous) restspace of the observer* at $s \in (a, b)$ is $T_{\dot{\gamma}(s)} \Sigma_{\gamma(s)}$. There are two natural metrics in this restspace. The first one is given by

$$g_{\dot{\gamma}(s)}|_\Sigma : T_{\dot{\gamma}(s)} \Sigma_{\gamma(s)} \times T_{\dot{\gamma}(s)} \Sigma_{\gamma(s)} \to \mathbb{R},$$

namely, the fundamental tensor restricted to $\Sigma$, which is a definite metric. The direction $\dot{\gamma}(s)$ is $g_{\dot{\gamma}(s)}$-orthogonal to $T_{\dot{\gamma}(s)}\Sigma$. As we have said above, this metric is the best approximation of $L$ with a scalar product in the direction of $\dot{\gamma}(s)$ as the restriction $g|_\Sigma$ is the second fundamental form of $\Sigma$ with respect to the opposite to the position vector (using the affine connection in $T_p M$).

The other metric is a Finsler metric with indicatrix $S_{\dot{\gamma}(s)}$, where $S_v = T_v \Sigma \cap \mathcal{C}_p$ (the set of velocities of light in the restspace of $v$). It is unclear which one is more suitable to measure spacelike distances, and indeed, the choice of metric could depend on the type of measure.

## 3   Anisotropic connections

### 3.1   Anisotropic tensor fields and their vertical derivatives

We will denote by $x = (x^1, \ldots, x^n)$ local coordinates on some open subset $U$ of $M$ and by $(x, y) = (x^1, \ldots, x^n, y^1, \ldots, y^n)$ the natural ones induced on $TU \subset TM$. Let us denote by $T^*M$ the cotangent bundle of the manifold $M$ and by $\bigotimes^{r)} TM \otimes \bigotimes^{s)} T^*M$ the classical vector bundle of tensors of type $(r, s)$ over $M$. Recall that an $(r, s)$-tensor field on $M$ is a smooth section of this bundle, and let $\mathcal{T}^r_s(M)$ be the space of all such tensor fields. We use the simplified notation for smooth functions and vector fields on $M$, resp., $\mathcal{T}^0_0(M) = \mathcal{F}(M)$ and $\mathcal{T}^1_0(M) = \mathfrak{X}(M)$.

Now, let $\pi^*_A(\bigotimes^{r)} TM \otimes \bigotimes^{s)} T^*M)$ be the bundle over $A$ pullbacked by the natural projection $\pi_A : A \to M$. A smooth section $T : A \to \pi^*_A(\bigotimes^{r)} TM \otimes \bigotimes^{s)} T^*M)$ of this bundle is called an *A-anisotropic tensor field* on $M$, and let $\mathcal{T}^r_s(M_A)$ the space of such fields (also by convention, $\mathcal{T}^0_0(M_A) = \mathcal{F}(A)$). In natural coordinates with summation in repeated indices,

$$T_v = T^{a_1, \ldots, a_r}_{b_1, \ldots, b_s}(x, y) \, \partial_{a_1}|_x \otimes \ldots \otimes \partial_{a_r}|_x \otimes dx^{b_1}|_x \otimes \ldots \otimes dx^{b_s}|_x \qquad \forall v \in A \cap TU; \ (1)$$

here, $(x, y)$ and $x$ are, resp., the coordinates of $v$ and $\pi_A(v)$, while the functions $T^{a_1, \ldots, a_r}_{b_1, \ldots, b_s}$ transform tensorially under changes of coordinates. Recall that, naturally, $\mathcal{T}^r_s(M_A)$ becomes a module over the ring $\mathcal{F}(A)$ and the tensor products and contractions induce further operations on sections, as in the case of usual tensor fields on $M$. In particular, $\mathcal{T}^1_0(M_A)$ and $\mathcal{T}^0_1(M_A)$ will be called resp. the sets of *anisotropic vector fields* and *1-forms* on $M$. We emphasize a particularity of the former: the elements $X \in \mathcal{T}^1_0(M_A)$ are sections of $\pi^*_A(TM) \to A$ and this is naturally isomorphic to the *vertical bundle* $\mathcal{V}A \to A$, where $\mathcal{V}_{(x,y)}A =$ $\mathrm{Span}_{\{i=1,\dots,n\}}\left\{\partial_{y^i}|_{(x,y)}\right\} \subset T_{(x,y)}A$. Thus, $X$ can be identified with a vertical vector field on $A$,

$$X_v = X^i(x,y)\partial_{x^i}|_x \in T_{\pi(v)}M \leftrightarrow X_v^{\mathcal{V}} = X^i(x,y)\partial_{y^i}|_{(x,y)} \in \mathcal{V}_v A. \qquad (2)$$

It is also worth pointing out that there is a natural inclusion

$$\mathcal{T}_s^r(M) \hookrightarrow \mathcal{T}_s^r(M_A), \qquad T \mapsto \tilde{T}, \qquad (3)$$

just putting the components of $\tilde{T}$ in (1) as independent of $y$ and equal to those of $T$. The tensor $\tilde{T}$ will be called *isotropic* and we will not distinguish between $T$ and $\tilde{T}$ when there is no possibility of confusion.

Notice that, at each $v \in A$, the fiber of the bundle $\pi_A^*(\overset{r)}{\bigotimes}TM \otimes \overset{s)}{\bigotimes}T^*M)$ becomes the space of all the $(r,s)$-tensors on $p = \pi_A(v)$. As this is a single vector space, the derivative of any curve of such tensors is well defined. Thus, given an anisotropic tensor $T \in \mathcal{T}_s^r(M_A)$ we can define its *vertical derivative* at $v \in A$ in any direction $w \in T_pM$, $p = \pi(v)$, as

$$(\dot{\partial}_w T)_v = \frac{d}{dt}T_{v+tw}|_{t=0},$$

which is again a tensor on $T_pM$ regarded as a fiber of $v$. As the map $w \mapsto (\dot{\partial}_w T)_v$ is linear, we can introduce naturally an $(r, s+1)$ tensor field as follows (we write directly the obvious expression in coordinates).

**Definition 3.** *Given an A-anisotropic tensor $T \in \mathcal{T}_s^r(M_A)$, its vertical derivative $\dot{\partial}T \in \mathcal{T}_{s+1}^r(M_A)$ is defined as*

$$(\dot{\partial}T)_v = \partial_{y^{b_{s+1}}}T_{b_1,\dots,b_s}^{a_1,\dots,a_r}(x,y)\,\partial_{a_1}|_x \otimes \dots \otimes \partial_{a_r}|_x \otimes \mathrm{d}x^{b_1}|_x \otimes \dots \otimes \mathrm{d}x_x^{b_s} \otimes \mathrm{d}x_x^{b_{s+1}}$$

*in any natural coordinates as in* (1).

### 3.2  Basic notion of anisotropic connection

As with other kinds of connections, anisotropic (linear) connections can be defined in different ways. We will introduce them by using a restricted domain which will be extended unequivocally later. Then we will refer to any of its characterizations also as an anisotropic (or, more properly, *A*-anisotropic) connection.

**Definition 4.** *An A-anisotropic connection is a map*

$$\nabla : \mathfrak{X}(M) \times \mathfrak{X}(M) \to \mathcal{T}_0^1(M_A), \qquad (X,Y) \mapsto \nabla_X Y,$$

*such that*

*(i)* $\nabla_X(Y + Z) = \nabla_X Y + \nabla_X Z$ *for any* $X, Y, Z \in \mathfrak{X}(M)$,
*(ii)* $\nabla_X(fY) = (X(f)Y) \circ \pi_A + (f \circ \pi_A)\nabla_X Y$ *for any* $f \in \mathcal{F}(M)$, $X, Y \in \mathfrak{X}(M)$,

*(iii)* $\nabla_{fX+hY}Z = (f \circ \pi_A)\nabla_X Z + (h \circ \pi_A)\nabla_Y Z$ *for any* $f, h \in \mathcal{F}(M)$ *and*
   $X, Y, Z \in \mathfrak{X}(M)$.

*We can eventually use the notation* $\nabla_X^v Y = (\nabla_X Y)_v$ *and, consistently,*

$$\nabla_X^V Y := (\nabla_X Y) \circ V \ (\in \mathfrak{X}(U)), \qquad \text{for any } A\text{-admissible vector field } V \in \mathfrak{X}(U).$$

*Moreover, we will say that the anisotropic connection is* homogeneous *(of degree zero), or* invariant by homotheties, *if for every* $v \in A$ *and* $\lambda > 0$, $(\nabla_X Y)_{\lambda v} = (\nabla_X Y)_v$ *(that is,* $\nabla_X Y = \nabla_X Y \circ h_\lambda$ *where* $h_\lambda : A \to A$ *is the homothety* $h_\lambda(v) = \lambda v$*).*

As in the case of affine connections, $\nabla$ has a local nature and, by using coordinates, we can express this in terms of the *Christoffel symbols of* $\nabla$, which are defined as the functions $\varGamma_{ij}^a : TU \cap A \to \mathbb{R}$ determined by

$$(\nabla_{\partial_i}\partial_j)_v \ (= \nabla_{\partial_i}^v \partial_j) = \varGamma_{ij}^a(v)\partial_a|_{\pi(v)}, \qquad \text{that is,} \qquad \nabla_{\partial_i}^V \partial_j = \varGamma_{ij}^a(V)\partial_a. \quad (4)$$

Clearly, the homogeneity of $\nabla$ is then equivalent to the 0-homogeneity of its Christoffel symbols, $\varGamma_{ij}^a(\lambda v) = \varGamma_{ij}^a(v)$, $\lambda > 0$. The following properties of these symbols are proven as in the standard case of affine connections.

**Proposition 1.** *(1) Under a change of coordinates* $(U, x) \rightsquigarrow (\bar{U}, \bar{x})$, *the Christoffel symbols* $\varGamma_{kl}^m$, $\bar{\varGamma}_{ij}^a$ *satisfy the cocycle transformation*

$$\bar{\varGamma}_{ij}^a(x, y) = \frac{\partial \bar{x}^a}{\partial x^m}(x)\left(\frac{\partial^2 x^m}{\partial \bar{x}^i \partial \bar{x}^j}(x) + \frac{\partial x^k}{\partial \bar{x}^i}(x)\,\frac{\partial x^l}{\partial \bar{x}^j}(x)\,\varGamma_{kl}^m(x, y)\right) \qquad (5)$$

*(2) Conversely, given any local choice of functions* $\varGamma_{ij}^k$ *for a coordinate atlas satisfying the cocycle transformation* (5), *there exists a unique anisotropic connection* $\nabla$ *whose Christoffel symbols are these functions. Moreover, if the functions are 0-homogeneous in* $y$, *then the produced* $\nabla$ *is homogeneous too.*

*(3) Any (classical, Koszul) affine connection on* $M$ *induces naturally an anisotropic one with Christoffel symbols independent of* $y$, *for any open conic domain* $A \subset TM$ *which naturally projects onto the whole* $M$.

*(4) Given an anisotropic connection* $\nabla$ *with Christoffel symbols* $\varGamma_{ij}^a$ *for each coordinates* $(U, x)$, *the choice of functions* $\varGamma_{ji}^a$ *for each* $(U, x)$ *yields a new connection* $\hat{\nabla}$, *and* $\nabla$ *is called* symmetric *if* $\nabla = \hat{\nabla}$.

### 3.3 Extended notion to take derivatives of anisotropic tensors

As a first observation, an anisotropic connection can be extended to an (anisotropic) covariant derivation of isotropic tensor fields in a natural way,

$$\nabla : \mathfrak{X}(M) \times \mathcal{T}_s^r(M) \to \mathcal{T}_s^r(M_A), \qquad (X, Y) \mapsto \nabla_X Y \qquad (6)$$

(also denoted $\nabla$ with no possiblity of confusion). This is done just defining

$$\nabla_X f = X(f) \circ \pi_A \qquad \text{(i.e. } \nabla_{\partial_i}f = \partial_i f \circ \pi_A\text{)} \qquad \forall f \in \mathcal{F}(M) = \mathcal{T}_0^0(M) \quad (7)$$

and using the usual rules of Leibniz and commutation with contractions. For example, if $\omega \in \mathcal{T}_1^0(M)$, then $\nabla_X \omega$ is determined by applying it to vector fields $\tilde{Z} \in \mathcal{T}_0^1(M_A)$ induced from $Z \in \mathcal{T}_0^1(M)$ as in (3), imposing

$$(\nabla_X \omega)(\tilde{Z}) = X(\omega(Z)) - \tilde{\omega}(\nabla_X Z)$$

(where $\tilde{\omega}_{(x,y)} = \omega_x \in T_x^* M$) or, in coordinates,

$$\nabla_{\partial_i}^V dx^k = -\Gamma_{ij}^k(V) \, dx^j.$$

Next, we will go beyond extending the anisotropic connection to an *(anisotropic) covariant derivation of A-anisotropic tensor fields* for any $r, s \geq 0$,

$$\nabla : \mathfrak{X}(M) \times \mathcal{T}_s^r(M_A) \to \mathcal{T}_s^r(M_A), \qquad (X, Y) \mapsto \nabla_X Y, \tag{8}$$

in a natural way (again, the same symbol $\nabla$ will be used without possibility of confusion, see the item *(i)* in Theorem 1 below). The key to get this new extension of $\nabla$ is to find a definition of $\nabla_X h$ when $h \in \mathcal{F}(A)$, that is, to find a natural extension of (7). This natural choice will be

$$(\nabla_X h)(v) = X_p(h \circ V) - \dot{\partial}_{(\nabla_{X_p} V)} h, \tag{9}$$

where $p = \pi_A(v)$ and $V \in \mathfrak{X}(M)$ satisfies $V_p = v$.

**Lemma 1.** *The definition of $\nabla_X h$ in (9) is independent of the choice of $V$. Moreover, if $h = f \circ \pi_A$ for some $f \in \mathcal{F}(A)$, then $\nabla_X h$ is equal to $\nabla_X f$ in (7).*

*Proof.* It is enough to check that the expression (9) written in coordinates is independent of $V$. Let $V = V^j \partial_j$, $X = X^i \partial_i \in \mathfrak{X}(M)$:

$$X(h \circ V) = X(h(x, V(x))) = X^i \left( \frac{\partial h}{\partial x^i}(x) + \frac{\partial h}{\partial y^j} \frac{\partial V^j}{\partial x^i} \right)(x, V(x)),$$
$$\dot{\partial}_{(\nabla_X V)} h = X(V^i) \, \dot{\partial}_i h + X^i V^j \Gamma_{ij}^k(x, V(x)) \, \dot{\partial}_k h,$$

the latter using that $\nabla_X V = X(V^i)\partial_i + X^i V^j \nabla_{\partial_i} \partial_j$. So, (9) reads at each $v$ of coordinates $(x, y)$:

$$(\nabla_X h)(x, y) = X^i(x) \frac{\partial h}{\partial x^i}(x, y) - X^i(x) \, y^j \, \Gamma_{ij}^k(x, y) \frac{\partial h}{\partial y^k} = X^i \left( \partial_i h - \Gamma_{ij}^k \, y^j \dot{\partial}_k h \right), \tag{10}$$

which is independent of the chosen $V$, as required.

Even though this lemma ensures the consistency of the definition of $\nabla_X h$, its meaning is not so evident. Algebraically, it ensures a sort of chain rule for $X(h \circ V)$. Anyway, we will give a further interpretation (see Remark 5). As a summary of this subsection, we obtain

**Theorem 1.** *For any anisotropic connection as in Definition 4 there exists a unique map as in (8) for each $r, s \geq 0$ satisfying, for any $T \in \mathcal{T}_s^r(M_A), T' \in \mathcal{T}_{s'}^{r'}(M_A)(r, r's, s' \geq 0)$,*

*(i)* $\nabla_X h$ *is given by* (10);

*(ii) Leibniz rule:* $\nabla_X(T \otimes T') = (\nabla_X T) \otimes T' + T \otimes (\nabla_X T')$;

*(iii) Commutativity with contractions:* $\nabla_X \left( C_j^i(T) \right) = C_j^i(\nabla_X T)$, *where* $C_j^i$ *denotes any possible contraction of* $T$ *in the slots* $i$ *(contravariant) and* $j$ *(covariant).*

Moreover, *if* $T$ *is written in coordinates as in* (1), *then*

$$T^{a_1,\ldots,a_r}_{b_1,\ldots,b_s|k} := (\nabla_k T)^{a_1,\ldots,a_r}_{b_1,\ldots,b_s}$$
$$= \partial_k T^{a_1,\ldots,a_r}_{b_1,\ldots,b_s} - \Gamma^i_{kj}\, y^j\, \dot\partial_i T^{a_1,\ldots,a_r}_{b_1,\ldots,b_s} + \textstyle\sum_{l=1}^r \Gamma^{a_l}_{kj_l} T^{a_1,\ldots,j_l,\ldots,a_r}_{b_1,\ldots,b_s} - \sum_{l=1}^s \Gamma^{i_l}_{kb_l} T^{a_1,\ldots,a_r}_{b_1,\ldots,i_l,\ldots b_s},$$

*where the coordinates depend on* $(x,y)$ *(i.e. on* $v \in TU \cap A$), $\partial_k = \partial_{x^k}, \dot\partial_k = \partial_{y^k}.$

The proof can be carried out by following the indications above. Anyway, full computations can be found in [10], where the following intrinsic version of the last displayed formula (regarding $T$ as a $\mathcal{F}(A)$-multilinear map) is proved [10, Theorem 11]: for any $v \in A$ and (local) extension $V \in \mathfrak{X}(U)$ of $v$,

$$
\begin{aligned}
(\nabla_X T)_v(\theta^1,\ldots,\theta^r,X_1,\ldots,X_s) = &\, X_{\pi(v)}(T_V(\theta^1,\ldots,\theta^r,X_1,\ldots,X_s)) \\
&- (\dot\partial T)_v(\theta^1,\ldots,\theta^r,X_1,\ldots,X_s,\nabla^V_X V), \\
&- \sum_{i=1}^r T_v(\theta^1,\ldots,\nabla_X\theta^i,\ldots,\theta^r,X_1,\ldots,X_s) \\
&- \sum_{j=1}^s T_v(\theta^1,\ldots,\theta^r,X_1,\ldots,\nabla_X X_j,\ldots,X_s),
\end{aligned}
$$
(11)

where $X,X_1,\ldots,X_s \in \mathfrak{X}(M)$ $(:= \mathcal{T}^1_0(M))$ and $\theta^1,\ldots,\theta^r \in \mathfrak{X}^*(M)$ $(:= \mathcal{T}^0_1(M))$.

*Remark 1.* Indeed, one can even extend $\nabla$ to a map

$$\nabla : \mathcal{T}^1_0(M_A) \times \mathcal{T}^r_s(M_A) \to \mathcal{T}^r_s(M_A), \qquad (X,Y) \mapsto \nabla_X Y,$$

just making it $\mathcal{F}(A)$-linear with respect to the first variable.

*Remark 2.* We have seen that the domain of an anisotropic connection can be extended from vector fields $X,Y \in \mathfrak{X}(M)$ to anisotropic vector fields in $\mathcal{T}^1_0(M_A)$. Additionally, multilinear maps over anisotropic tensor fields valued on anisotropic vector fields can be regarded as anisotropic tensor fields.

A relevant example appears when two anisotropic connections $\bar\nabla, \nabla$ are considered. Its difference $Q = \bar\nabla - \nabla$ is naturally an $\mathcal{F}(M)$-multilinear map

$$\mathfrak{X}(M) \times \mathfrak{X}(M) \longrightarrow \mathcal{T}^1_0(M_A),$$

which extends unequivocally by $\mathcal{F}(A)$-multilinearity to an anisotropic tensor field $Q \in \mathcal{T}^1_2(M_A)$ (recall $\mathfrak{X}(M) = \mathcal{T}^1_0(M) \hookrightarrow \mathcal{T}^1_0(M_A)$ as in (3)). Moreover, $Q$ can also be regarded as a multilinear map

$$\mathcal{T}^1_0(M_A) \times \mathcal{T}^1_0(M_A) \to \mathcal{T}^1_0(M_A), \text{ or } \mathcal{T}^1_0(M_A) \times \mathcal{T}^1_0(M_A) \times \mathcal{T}^0_1(M_A) \to \mathcal{F}(A).$$

In particular, its components $Q^k_{ij}$ share the cocycle transformation of a tensor field (i.e. the rule (5) without the term in the second derivatives $\partial^2 x^m/\partial\bar x^i\partial\bar x^j$).

Applying the previous discussion to the case of the connection $\hat{\nabla}$ obtained from $\nabla$ with Christoffel symbols $\hat{\Gamma}_{ji}^{k} = \Gamma_{ij}^{k}$ (see Proposition 1), the following definition becomes consistent.

**Definition 5.** *The* torsion *of an anisotropic connection $\nabla$ is the tensor* Tor $\in \mathcal{T}_2^1(M_A)$ *defined in coordinates by:*

$$\text{Tor}_{ij}^k = \Gamma_{ij}^k - \Gamma_{ji}^k.$$

The following consequence is also straightforward now.

**Corollary 1.** *The space of all the anisotropic connections on $M$ has a structure of affine space with associated vector space $\mathcal{T}_2^1(M_A)$.*

*Moreover, the homogeneous anisotropic connections become an affine subspace with associated vector space the subspace of $\mathcal{T}_2^1(M_A)$ composed by its 0-homogeneous tensors (i.e., satisfying $Q_{ij}^k(x, \lambda y) = Q_{ij}^k(x, y)$).*

## 4   Anisotropic vs nonlinear connections

Any anisotropic connection induces a nonlinear connection on $\pi_A : A \to M$. Let us start recalling the framework of the latter adapted to our case.

As $A \subset TM$ then $TA \subset T(TM)$ and the coordinates $(x, y)$ on $TU$ induce naturally coordinates $(x, y, \dot{x}, \dot{y})$ for some neighborhood in $TA$. The *vertical bundle* $\mathcal{V}A \subset TA$ is the subbundle of $TA \to A$ composed by the elements with $\dot{x} = 0$, so that $(x, y, \dot{y})$ become natural coordinates. Clearly, $\mathcal{V}A$ is naturally identifiable to the pullbacked bundle of vectors $\pi_A^*(TM)$ and, consistently, the sections of the bundle $\mathcal{V}A \to A$ are naturally identifiable to the space of anisotropic vector fields $\mathcal{T}_0^1(M_A)$.

### 4.1   Standard setting for nonlinear connections.

There are several ways to define a nonlinear connection on $A \to M$. One of them is to provide a vector bundle homomorphism $\nu : TA \to \mathcal{V}A$ such that $\nu|_{\mathcal{V}A}$ is the identity. Then, the distribution $\mathcal{H}A := \ker(\nu)$ characterizes $\nu$ and gives a decomposition $TA = \mathcal{H}A \oplus \mathcal{V}A$, which can also be used as an alternative definition of the nonlinear connection. In coordinates:

$$\nu \left( \dot{x}^i \frac{\partial}{\partial x^i} \bigg|_{(x,y)} + \dot{y}^a \frac{\partial}{\partial y^a} \bigg|_{(x,y)} \right) = \left( \dot{y}^a + N_i^a(x,y)\dot{x}^i \right) \frac{\partial}{\partial y^a} \bigg|_{(x,y)}$$

$$\mathcal{H}_{(x,y)}A = \text{Span}_{\{i=1,\ldots,n\}} \left\{ \frac{\delta}{\delta x^i} \bigg|_{(x,y)} := \frac{\partial}{\partial x^i} \bigg|_{(x,y)} - N_i^a(x,y) \frac{\partial}{\partial y^a} \bigg|_{(x,y)} \right\},$$

where the $N_i^a$ 's are defined on $A \cap TU$ and, taking a new chart $\bar{U}, (\bar{x}, \bar{y})$:

$$\bar{N}_j^a(\bar{x}, \bar{y}) = \frac{\partial x^i}{\partial \bar{x}^j} \left( -\frac{\partial^2 \bar{x}^a}{\partial x^i \partial x^l} y^l + \frac{\partial \bar{x}^a}{\partial x^b} N_i^b(x,y) \right) \qquad (12)$$

on $U \cap \bar{U}$.

*Remark 3.* (a) Conversely, any covering of charts of $M$ endowed with a set of functions satisfying the cocycle transformation (12) determines unequivocally a nonlinear connection of $A \to M$. Incidentally, we recover a standard fact from the theory of fibered manifolds: a nonlinear (or *Ehresmann*) connection on $A \to M$ is the same as a section of the 1-jet bundle $\mathbf{J}^1 A \to A$; see [22, §17] for instance. If, by means of such a section, for the selected 1-jet at $v \in A$ one puts

$$N_i^a(v) = -\frac{\partial V^a}{\partial x^i}(\pi(v))$$

(where $V$ is a local extension of $v$ that determines the jet), then it is straightforward to see that these $N_i^a$'s satisfy (12).

(b) It makes sense to assume that the nonlinear connection on $A \to M$ is (positive) *homogeneous* in the sense that the distribution $\mathcal{H}A$ is invariant by homotheties $h_\lambda : A \to A$, that is, if $Th_\lambda : TA \to TA$ is the tangent map (differential) of $h_\lambda$, then for each $\lambda > 0$, $(x, y) \equiv v \in A$:

$$(Th_\lambda)_v(\mathcal{H}_v) = \mathcal{H}_{\lambda v} \qquad \text{or, equivalently,} \qquad N_i^a(x, \lambda y) = \lambda N_i^a(x, y). \tag{13}$$

For the latter, as an integral curve $(x(t), y(t))$ of the horizontal vector $\delta/\delta_{x^i}$ satisfies $dy^a/dx^i = -N_i^a(x, y)$, the 1-homogeneity of the functions $N_i^a$ characterizes when $(x(t), \lambda y(t))$ $(= h_\lambda(x(t), y(t)))$ is also an integral curve. In this case, the relations

$$\delta_{x^i}|_{\lambda v} = (Th_\lambda)_v (\delta_{x^i}|_v), \qquad \partial_{y^i}|_{\lambda v} = \lambda^{-1} (Th_\lambda)_v (\partial_{y^i}|_v)$$

(0-homogeneity of $\delta_{x^i}$ and $(-1)$-homogeneity of $\partial_{y^i}$, see [4, §1.5]) are also satisfied,

$$N_i^a(x, y) = \partial_{y^j} N_i^a(x, y) y^j \tag{14}$$

(by Euler theorem) and $\partial_{y^j} N_i^a(x, y)$ is 0-homogeneous in $y$ (i.e., invariant under $h_\lambda : (x, y) \mapsto (x, \lambda y)$).

**The case of linear connections.** As usual, the name "nonlinear" must be understood in the sense of "non-necessarily linear" for a connection. To introduce the linear case, recall first the following elementary technical result.

**Lemma 2.** *Assume that the functions $N_i^a : A \cap TU \to \mathbb{R}$ are 1-homogeneous in $y$ and can be smoothly extended to 0. Then there exist smooth functions $\Gamma_{ij}^k : U \to \mathbb{R}$ such that*

$$N_i^a(x, y) = \Gamma_{ij}^a(x) y^j. \tag{15}$$

*In particular, $N_{i\cdot j}^a(x) := \partial_{y^j} N_i^a(x, y) = \Gamma_{ij}^a(x)$ and each $N_i^a$ can be naturally extended to $TU$.*

*Proof.* Just recall that at each $x_0 \in U$ the function $f(y) = N_i^a(x_0, y)$ is linear, as $df|_0(v) = \lim_{\lambda \searrow 0} f(\lambda v)/\lambda = f(v)$ for each $v \equiv (x_0, y) \in A_{x_0}$ by homogeneity. $\square$

**Definition 6.** *We say that a nonlinear connection on $A \to M$ is* linear *when its horizontal distribution can be smoothly extended to the zero section (i.e., when its coefficients $N_i^a$'s satisfy* (15) *in every coordinate chart).*

*Remark 4.* (1) For linear connections, the dependence of $\mathcal{H}_{(x,y)}$ in $y$ is usually dropped.

(2) As suggested at the end of Lemma 2, any linear connection on $A \to M$ can be extended unequivocally to the tangent bundle. So, linear connections are usually defined on vector bundles.

Most of the framework of linear connections explained here for $A \to M$ is extended to more general bundles in a standard way and will be used without further mention in §6.

### 4.2 Projection of anisotropic connections onto nonlinear ones

Consider an anisotropic connection $\nabla$ as in the original Definition 4. It can be proved that $\nabla$ induces a horizontal distribution and, then, a nonlinear connection by using an intrinsic approach [10, §3.1]. However, we would like to emphasize the following direct relation with the cocycle trasformation associated with the Christoffel symbols $\Gamma_{ij}^a(x,y)$ in (4).

**Theorem 2.** *(1) Any anisotropic connection $\nabla$ defines canonically a nonlinear connection $\nu^\nabla$ whose components $N_i^a$ in any chart are given as[3]:*

$$N_i^a(x,y) = \Gamma_{ij}^a(x,y)\, y^j. \tag{16}$$

*Moreover, if $\nabla$ is homogeneous then $\nu^\nabla$ is homogeneous.*

*(2) Any nonlinear connection $\nu$ defines canonically an anisotropic connection $\nabla^\nu$ whose symbols $N_{ij}^a$ in any chart are given by:*

$$\Gamma_{ij}^a(x,y) = \dot\partial_j N_i^a(x,y) \,(= N_{i,j}^a(x,y)).$$

*Moreover, if $\nu$ is homogeneous then $\nabla^\nu$ is homogeneous and $\nu^{(\nabla^\nu)} = \nu$. In particular, the map $\nabla \mapsto \nu^\nabla$ is onto when restricted to the sets of homogeneous anisotropic and nonlinear connections.*

*(3) For any homogeneous $\nu$, the set $\{\nabla \text{ anisotropic connection} : \nu = \nu^\nabla\}$ is:*

$$\{\nabla^\nu + Q : \quad Q \in \mathcal{T}_2^1(M_A) \text{ and } Q_{ij}^k y^j = 0\}.$$

*(4) When restricted to (classical) affine and linear connections, the map $\nabla \mapsto \nu^\nabla$ is well defined and bijective.*

---

[3] Indeed, there would also be a second one $\Gamma_{ij}^a(x,y)\, y^i$ which differs in the sign of its torsion (recall Remark 2).

*Proof.* (1) Notice that the functions $N_i^a$ satisfy the cocycle transformation (12). Then, the 0-homogeneity in $y$ of $\Gamma_{ij}^a(x,y)$ gives the 1-homogeneity of $N_i^a(x,y)$.

(2) The cocycle (12) for $N_i^a(x,y)$ implies the cocycle (5) for $\Gamma_{ij}^a$. Then, the homogeneity of $\nu$ implies (14) and thus (16).

(3) Straightforward from part (2) and (16).

(4) For such a $\nabla$, $\Gamma_{ij}^k(x,y) \equiv \Gamma_{ij}^k(x)$ and $\nu^\nabla$ is linear. Given a second $\bar\nabla$ with $\nu^{\bar\nabla} = \nu^\nabla$, the difference $Q = \bar\nabla - \nabla$ is also independent of $y$ and $Q_{ij}^k(x)y^j = 0$, which implies $Q = 0$ by taking derivatives with respect to each $\partial_{y^i}$.    $\square$

*Remark 5.* Finally, we can give the promised interpretation of the definition of $\nabla_X h$ in (9), (10). Indeed, any $X_p = X^i(x)\partial_i|_x \in T_pM$, $p \equiv x$, gives a horizontal lift $X_v^{\mathcal H} = X^i(x)\,\delta/\delta x^i|_{(x,y)} \in \mathcal{H}_v A$ $(\subset T_v A)$ for any $v \equiv (x,y)$ with $\pi_A(v) = p$. So, substituting the expressions in the formulas of $\mathcal{H}_{(x,y)}A$ and (16):

$$X_v^H(h) = X^i(\pi(v))\left(\partial_i h - \Gamma_{ij}^k\, y^j \dot\partial_k h\right)(v) = (\nabla_X h)_v,$$

the last equality in agreement with (10).

## 5    Anisotropic vs linear connections

Next, our goal will be to identify anisotropic connections with a class of linear connections on the vector bundle $\mathcal{V}A \to A$.

As shown at the beginning of §4, $\mathcal{V}A$ $(\subset TA)$ admits as natural coordinates $(x,y,\dot y)$, which here will be relabelled $(x,y,z)$. Thus, we can write $z^a \partial_{y^a}|_{(x,y)} \in \mathcal{V}_{(x,y)}A$. The tangent space $T(\mathcal{V}A)$ includes a new vertical space $\mathcal{V}(\mathcal{V}A) \to \mathcal{V}A$ composed by the vectors of the type $\dot z^a \partial_{z^a}|_{(x,y,z)}$ .

### 5.1    Linear connections on $\mathcal{V}A \to A$

In order to define an *Ehresmann connection* on $\mathcal{V}A \to A$, we have to provide a smooth horizontal decomposition $T_{(x,y,z)}(\mathcal{V}A) = \mathcal{H}_{(x,y,z)} \oplus \mathcal{V}_{(x,y,z)}$ (just as is done in the case of $A \to M$, for which the Ehresmann connection is usually called *nonlinear connection*). Moreover, any positive homothety on $A$, $h_\lambda \colon (x,y) \mapsto (x,\lambda y)$, induces a natural morphism:

$$h_{\lambda_*} \colon \mathcal{V}A \to \mathcal{V}A, \quad (x,y,z) \mapsto (x, \lambda y, \lambda z)$$

(the restriction of $Th_\lambda$). As a natural extension of Remark 3(b), the new horizontal distribution, (and then the Ehresmann connection itself) is called *invariant by homotheties* when it is preserved by the tangent map of $h_{\lambda_*}$, that is,

$$(Th_{\lambda_*})_{(x,y,z)}\left(\mathcal{H}_{(x,y,z)}\right) = \mathcal{H}_{(x,\lambda y, \lambda z)}. \tag{17}$$

In what follows, we will focus on the particular case of a *linear connection* $\nu^*$ on $\mathcal{V}A \to A$. In this case, the horizontal decomposition is naturally independent of $z$ (extending Remark 4), so one can write:

$$T_{(x,y,z)}(\mathcal{V}A) = \mathcal{H}_{(x,y)} \oplus \mathcal{V}_{(x,y)}. \tag{18}$$

Consistently, in the case that $\nu^*$ is invariant by homotheties, the dependence of $\mathcal{H}$ on $z$ is dropped in (17). A standard result of the theory of linear connections says that $\nu^*$ is characterized by its *Koszul derivative* [4] $\nabla^*$, that is, the way the sections of $\mathcal{V}A \to A$ are covariantly differentiated (Theorems 2.52 and 2.58 on [21]). As these sections are identified with the anisotropic vector fields $\mathcal{T}_0^1(M_A)$ (recall the vertical isomorphism (2)), $\nabla^*$ becomes a map

$$\nabla^* : \mathfrak{X}(A) \times \mathcal{T}_0^1(M_A) \to \mathcal{T}_0^1(M_A), \qquad (W, Z) \mapsto \nabla_W^* Z \qquad (19)$$

(the reader may appreciate the similarities and differences between (19) and the anisotropic covariant derivative (8), the latter with $(r, s) = (1, 0)$). Moreover, it is straightforward but tedious to prove (from the definition of $\nabla^*$ in terms of the corresponding horizontal decomposition) that the invariance by homotheties of $\nabla^*$ is characterized the following way. If the section $Z^{\mathcal{V}} : A \to \mathcal{V}A$ is 0-homogeneous (meaning $Z^{\mathcal{V}} \circ h_\lambda = \lambda^{-1} h_{\lambda_*} \circ Z^{\mathcal{V}}$ as in [4, §1.5]), then

$$\nabla^*_{(Th_\lambda)_{(x,y)}(W_{(x,y)})} Z^{\mathcal{V}} = \lambda^{-1} (Th_\lambda)_{(x,y)} \left( \nabla^*_{W_{(x,y)}} Z^{\mathcal{V}} \right). \qquad (20)$$

To specify $\nabla^*$ by means of its Christoffel symbols, one has to choose a basis for $\mathfrak{X}(A)$ and another one for $\mathcal{T}_0^1(M_A)$. A possible choice would be the one associated to coordinates, namely $\{\partial_i|_{(x,y)}, \dot{\partial}_i|_{(x,y)}\} = \{\partial_{x^i}|_{(x,y)}, \partial_{y^i}|_{(x,y)}\}$ for the former and [5] $\{\partial_j|_x\} \equiv \{\dot{\partial}_j|_{(x,y)}\}$ for the latter. However, in case we have a prescribed nonlinear connection $\overset{o}{\nu}$ on $A \to M$, a more convenient choice than $\{\partial_i, \dot{\partial}_i\}$ may be $\{\delta_i, \dot{\partial}_i\}$:

$$\delta_i|_{(x,y)} = \frac{\delta}{\delta x^i}\bigg|_{(x,y)} = \partial_{x^i}|_{(x,y)} - \overset{o}{N}{}^a_i(x,y)\partial_{y^a}|_{(x,y)}, \qquad \dot{\partial}_j|_{(x,y)} = \partial_{y^j}|_{(x,y)}.$$

This happens in the pseudo-Finsler case, where $\overset{o}{\nu}$ is provided by the geodesic spray and, thus, is homogeneous. This last property of $\overset{o}{\nu}$ will be assumed for the sake of simplicity, even though actually it will be important only when the homogeneity of the Christoffel symbols is involved. From the homogeneity of $\overset{o}{\nu}$, $\delta_i$ is 1-homogeneous, namely $\delta_i|_{(x,\lambda y)} = (Th_\lambda)_{(x,y)} (\delta_i|_{(x,y)})$, while $\dot{\partial}_i$ is 0-homogeneous, namely $\dot{\partial}_i|_{(x,\lambda y)} = \lambda^{-1} (Th_\lambda)_{(x,y)} (\dot{\partial}_i|_{(x,y)})$.

**Definition 7.** *The horizontal and vertical Christoffel symbols of $\nabla^*$ with respect to a prescribed homogeneous nonlinear connection $\overset{o}{\nu}$ in the coordinates $(U, x)$ are, resp., the functions $H^a_{ij}$, $V^a_{ij}$ on $A \cap TU$ determined by*

$$H^a_{ij}(x,y)\,\partial_a|_x = \nabla^*_{\delta_i|_{(x,y)}} \partial_j (\equiv \nabla^*_{\delta_i|_{(x,y)}} \dot{\partial}_j = H^a_{ij}(x,y)\dot{\partial}_a|_{(x,y)}),$$

$$V^a_{ij}(x,y)\,\partial_a|_x = \nabla^*_{\dot{\partial}_i|_{(x,y)}} \partial_j (\equiv \nabla^*_{\dot{\partial}_i|_{(x,y)}} \dot{\partial}_j = V^a_{ij}(x,y)\dot{\partial}_a|_{(x,y)}).$$

---

[4] Such a result is similar to Theorem 2, item (4).

[5] Keep in mind that this identification also corresponds exactly to the vertical isomorphism (2).

**Proposition 2.** *(1) The Christoffel symbols $H_{ij}^a$, $V_{ij}^a$ of $\nabla^*$ with respect to $\overset{o}{\nu}$ satisfy:*

*(a) The cocycle for $H_{ij}^a$ (resp., $V_{ij}^a$) under a change of coordinates coincides with the one for the Christoffel symbols $\Gamma_{ij}^a$ of an A-anisotropic connection (5) (resp., the one of an A-anisotropic $(1,2)$ tensor). In particular, if all the $V_{ij}^a$'s vanish for some coordinates on $U$, then they vanish for any coordinates therein.*

*(b) In case that the linear connection $\nabla^*$ is invariant by homotheties, then:*

*(b1) $H_{ij}^a(x, \lambda y) = H_{ij}^a(x, y)$ (0-homogeneous in $y$), and*

*(b2) $V_{ij}^a(x, \lambda y) = \lambda^{-1} V_{ij}^a(x, y)$ ($(-1)$-homogeneous in $y$).*

*(2) Conversely, once a homogeneous nonlinear connection $\overset{o}{\nu}$ is prescribed, any local choice of functions $H_{ij}^a$, $V_{ij}^a$ satisfying (a) for a coordinate atlas determines a unique linear connection $\nabla^*$, whose Christoffel symbols with respect to $\overset{o}{\nu}$ in that atlas coincide with the original $H_{ij}^a$, $V_{ij}^a$. Moreover, if the functions $H_{ij}^a$, $V_{ij}^a$ are chosen as $0$ and $(-1)$-homogeneous in $y$ (consistently with (b) above), then the produced $\nabla^*$ is invariant by homotheties.*

*Proof.* (1) The cocycles of $H_{ij}^a$ and $V_{ij}^a$ can be checked from their definitions together with the transformation laws

$$\bar\delta_i = \frac{\partial x^k}{\partial \bar x^i}\delta_k, \qquad \bar{\dot\partial}_i = \frac{\partial x^k}{\partial \bar x^i}\dot\partial_k, \qquad \partial_j = \frac{\partial \bar x^l}{\partial x^j}\bar\partial_l$$

(recall that these last $x^j$'s are the ones on $M$ and not those of $TM$). Then, using only the definition of $H_{ij}^a$ and $V_{ij}^a$, the 1-homogeneity of the $\delta_i$'s, and the 0-homogeneity of the $\dot\partial_i$'s, the following identities hold true:

$$H_{ij}^a(x, \lambda y)\,\dot\partial_a\Big|_{(x,\lambda y)} = \nabla^*_{\delta_i|_{(x,\lambda y)}}\dot\partial_j = \nabla^*_{(Th_\lambda)_{(x,y)}(\delta_i|_{(x,y)})}\dot\partial_j, \tag{21}$$

$$H_{ij}^a(x, y)\,\dot\partial_a\Big|_{(x,y)} = \nabla^*_{\delta_i|_{(x,y)}}\dot\partial_j$$
$$= H_{ij}^a(x,y)\lambda\,(Th_\lambda)^{-1}_{(x,y)}\,(\dot\partial_a\Big|_{(x,\lambda y)}) = (Th_\lambda)^{-1}_{(x,y)}\,(\lambda H_{ij}^a(x,y)\,\dot\partial_a\Big|_{(x,\lambda y)}). \tag{22}$$

Now, in case that $\nabla^*$ is invariant by homotheties, one can use (20) with $Z^{\mathcal V} = \dot\partial_j$. From (20), (21) and (22) follows that

$$H_{ij}^a(x, \lambda y)\,\dot\partial_a\Big|_{(x,\lambda y)} = \lambda^{-1}\,(Th_\lambda)_{(x,y)}\,(\nabla^*_{\delta_i|_{(x,y)}}\dot\partial_j)$$
$$= \lambda^{-1}\,(Th_\lambda)_{(x,y)}\,(H_{ij}^a(x,y)\,\dot\partial_a\Big|_{(x,y)})$$
$$= \lambda^{-1}\,(Th_\lambda)_{(x,y)}\,\{(Th_\lambda)^{-1}_{(x,y)}\,(\lambda H_{ij}^a(x,y)\,\dot\partial_a\Big|_{(x,\lambda y)})\}$$
$$= H_{ij}^a(x,y)\,\dot\partial_a\Big|_{(x,\lambda y)},$$

thus proving (b1). An analogous reasoning proves (b2).

(2) Knowing the cocycles that the Christoffel symbols of such a $\nabla^*$ should satisfy, it is possible to define $\nabla^*$ by $H_{ij}^a$, $V_{ij}^a$ on each coordinate chart and, as usual, assert that the local definitions patch together to form a global linear connection on $\mathcal{V}A \to A$. Moreover, (21) and (22) are still valid for this $\nabla^*$, and so are the analogous identities for the $V_{ij}^a$'s. So, if the $H_{ij}^a$'s are 0-homogeneous and the $V_{ij}^a$'s are $(-1)$-homogeneous, then one can use those identities to show that

$$\nabla^*_{(Th_\lambda)_{(x,y)}(\delta_{i|_{(x,y)}})}\dot{\partial}_j = \lambda^{-1}\,(Th_\lambda)_{(x,y)}\,(\nabla^*_{\delta_{i|_{(x,y)}}}\dot{\partial}_j), \tag{23}$$

$$\nabla^*_{(Th_\lambda)_{(x,y)}(\dot{\partial}_{i|_{(x,y)}})}\dot{\partial}_j = \lambda^{-1}\,(Th_\lambda)_{(x,y)}\,(\nabla^*_{\dot{\partial}_{i|_{(x,y)}}}\dot{\partial}_j). \tag{24}$$

By using that $\{\delta_i, \dot{\partial}_i\}$ is a basis for $\mathfrak{X}(A)$, $\{\dot{\partial}_i\}$ is a basis for the vertical vector fields, and the two expressions $\nabla^*_{(Th_\lambda)_{(x,y)}(W_{(x,y)})}Z^{\mathcal{V}}$ and $\lambda^{-1}\,(Th_\lambda)_{(x,y)}\,(\nabla^*_{W_{(x,y)}}Z^{\mathcal{V}})$ are linear in $W$ and Leibnizian in $Z^{\mathcal{V}}$, the identities (23) and (24) prove that $\nabla^*$ satisfies (20), hence that $\nabla^*$ is invariant by homotheties. □

### 5.2   Anisotropic connections as vertically trivial linear connections.

For a linear connection $\nabla^*$ on $\mathcal{V}A \to A$, the vanishing of all the $V_{ij}^a$'s involves only vertical derivatives and, so, it is an intrinsic property (independent also of $\overset{o}{\nu}$ in Proposition 2 item (1)(a)). This makes the following consistent.

**Definition 8.** *Let $\nabla^*$ be a linear connection on $\mathcal{V}A \to A$. We say that $\nabla^*$ is* vertically trivial *if $V_{ij}^a = 0$ everywhere.*

*Remark 6.* From Proposition 2, it is clear that any homogeneous nonlinear connection $\overset{o}{\nu}$ induces a projection of the set of all the $\nabla^*$'s onto the vertically trivial ones, namely $(H_{ij}^a, V_{ij}^a) \mapsto (H_{ij}^a, 0)$.

**Theorem 3.** *Let $\overset{o}{\nu}$ be a homogeneous nonlinear connection on $A \to M$. The map between the sets of vertically trivial and $A$-anisotropic connections, defined in natural coordinates as*

$$\{\text{vertically trivial conn. } \nabla^* \text{ on } \mathcal{V}A \to A\} \longrightarrow \{A\text{-anisotropic conn. } \nabla\},$$
$$(H_{ij}^a, V_{ij}^a = 0) \qquad\qquad \mapsto \qquad \Gamma_{ij}^a = H_{ij}^a,$$

*is well defined, bijective and it maps (also bijectively) the homothetically invariant $\nabla^*$'s into the homogeneous (homothetically invariant) $\nabla$'s.*

*Moreover, this map does not depend on the chosen $\overset{o}{\nu}$. So, there exists a natural identification between vertically trivial connections and anisotropic connections, preserving the homothetic invariance property.*

*Proof.* The first part is straightforward from Proposition 2. For the independence of $\overset{o}{\nu}$, recall that when chosing a second $\nu'$, the differences $\delta_i' - \delta_i$ are vertical and, thus, $\nabla^*_{\delta_i' - \delta_i}\partial_j = 0$, as $\nabla^*$ is vertically trivial. For the last assertion, recall that $A \to M$ always admits a homogeneous nonlinear connection (for example, the associated with any pseudo-Finsler metric on $A$ and, in particular, with any Riemannian metric on $M$). □

## 6   Anisotropic vs Finsler connections

When a pseudo-Finsler metric $L$ is given on $A$, the standard approach focuses on two geometric structures. The first one is its geodesic spray on $A$ and, thus, its associated homogeneous nonlinear connection on $A \to M$; these are unequivocally constructed from $L$. The second one is a linear connection on $\mathcal{V}A \to A$ invariant by homotheties. However, a priori there is not a unique canonical choice for the latter. Let us explain briefly the interplay of anisotropic connections with these two structures.

Recall that, as $L$ is 2-homogeneous, the nonlinear, linear and anisotropic connections to be considered here will be homogeneous (invariant by homotheties). In particular, all the the conclusions of Theorem 2 will be applicable and, for example, the torsion of a nonlinear connection $\nu$ can be defined as the torsion of the corresponding anisotropic connection $\nabla^\nu$.

### 6.1   The geodesic spray.

**Definition 9.** *A* spray *on $A$ is a vector field $G$ there (section $A \to TA$) satisfying: (a) $G$ is a* second order equation, *that is, it can be written as[6]:*

$$G_{(x,y)} = y^i \left.\frac{\partial}{\partial x^i}\right|_{(x,y)} - 2G^a(x,y) \left.\frac{\partial}{\partial y^a}\right|_{(x,y)},$$

*and (b)* 2-homogeneity, *i.e., $G^a(x,\lambda y) = \lambda^2 G^a(x,y)$ for[7] $\lambda > 0$.*

The spray associated with a pseudo-Finsler $L$ is its *geodesic spray*, which comes from the Euler-Lagrange equation for the geodesics (regarded as critical points of the $L$-energy functional).

We summarize the relations between sprays and homogeneous nonlinear connections (analogous to Theorem 2) following [20, §3.4].

**Proposition 3.** *(1) A homogeneous nonlinear connection $\nu$ defines a natural spray $G = \mathbb{C}^{\mathcal{H}}$ (the $\nu$-horizontal lift of the Liouville $\mathbb{C}$), in coordinates*

$$G^a(x,y) = \frac{1}{2}N_i^a(x,y)y^i.$$

*The integral curves of $G$ are called* geodesics *of $\nu$.*

*(2) A spray $G$ yields a natural homogeneous nonlinear connection $\overset{o}{\nu}$ on the fibered manifold $A \to M$ with components*

$$\overset{o}{N_i^a}(x,y) = \frac{\partial G^a}{\partial y^i}(x,y) = G_{\cdot i}^a(x,y).$$

---

[6] More intrinsically, $T\pi_A \circ G$ is the identity in $A$. Recall that its integral curves $(x(t),y(t))$, usually called *geodesics of the spray*, satisfy $dx^i/dt = y^i, dy^a/dt = -2G^a(x,y)$ and, thus, $d^2x^a/dt^2 + 2G^a(x,dx/dt) = 0$.

[7] More intrinsically, $[\mathbb{C},G] = G$, where $\mathbb{C}_v = y^i\partial_{y^i}|_{(x,y)}$ is the canonical (Liouville) vector field on $A$.

Then, $G = \mathbb{C}^{\overset{\circ}{\mathcal{H}}}$ ($\overset{\circ}{\nu}$-horizontal lift of $\mathbb{C}$) and $\overset{\circ}{\nu}$ is torsion-free.

(3) The geodesics of a homogeneous nonlinear connection $\nu$ are the integral curves of a spray $G$ iff $G = \mathbb{C}^{\mathcal{H}}$.

(4) Given a homogeneous nonlinear connection $\nu$, consider the natural spray $G = \mathbb{C}^{\mathcal{H}}$ and the nonlinear connection $\overset{\circ}{\nu}$ associated with this $G$. Then the difference $\nu - \overset{\circ}{\nu}$ is in $\mathcal{T}_1^1(M_A)$ with components

$$N_i^a(x,y) - G_{\cdot i}^a(x,y) = \frac{1}{2}Tor_{ij}^a(x,y)y^j,$$

where $Tor$ is the torsion of $\nabla^\nu$ (see part (2) of Theorem 2). Moreover, if this difference vanishes, then actually $Tor$ vanishes.

(5) Any homogeneous nonlinear connection is determined by its geodesics and torsion.

*Proof.* (1) These $G^a$'s satisfy the cocycle transformation required to form a second order equation and their 2-homogeneity comes from the 1-homogeneity of $N_i^a$.

(2) From the cocycle of a second order equation, these $\overset{\circ}{N_i^a}$'s satisfy (12) and their 1-homogeneity comes from the 2-homogeneity of the $G^a$'s. Then $G = \mathbb{C}^{\overset{\circ}{\mathcal{H}}}$ is nothing but the Euler relation for the $G^a$'s, while the torsion of $\overset{\circ}{\nu}$ is given by $G_{\cdot i \cdot j}^a - G_{\cdot j \cdot i}^a = G_{\cdot i \cdot j}^a - G_{\cdot i \cdot j}^a = 0$.

(3) Recall that the $N$-geodesic equation is

$$\frac{dy^a}{dt} + N_i^a(x(t),y(t))y^i(t) = 0,$$

so these are the integral curves of $G$ if and only if $N_i^a(x,y)y^i = 2G^a(x,y)$.

(4) This is a straightforward computation taking into account that $G^a = N_i^a y^i/2$ and $Tor_{ij}^a = N_{i \cdot j}^a - N_{j \cdot i}^a$. Thus, if $Tor_{ij}^a y^j = 0$, then $\nu = \overset{\circ}{\nu}$ and its torsion vanishes due to (3).

(5) This follows directly from (4).     □

*Remark 7.* As a consequence of Theorem 2, the $\overset{\circ}{\nu}$ associated with $G$ can be always obtained from a canonical homogeneous anisotropic connection $\nabla^{\overset{\circ}{\nu}}$ (item (2) of the theorem). This $\nabla^{\overset{\circ}{\nu}}$ is the so-called *Berwald anisotropic connection*. The remaining anisotropic connections that yield $G$ would be controlled by an anisotropic tensor $Q$ satisfying $Q_{ij}^a y^i y^j = 0$.

The components of the geodesic spray of a pseudo-Finsler metric are $G^a = \gamma_{ij}^a y^i y^j$, where

$$\gamma_{ij}^a = \frac{1}{2}g^{ak}\left(\frac{\partial g_{ki}}{\partial x^j} + \frac{\partial g_{kj}}{\partial x^i} - \frac{\partial g_{ij}}{\partial x^k}\right)$$

are the so-called *formal Christoffel symbols*, (see for example [24, formula (4.30)]). From them, the $\overset{\circ}{\nu}$ of Prop. 3 (2) is the connection given by

$$\overset{\circ}{N_i^a} = \gamma_{ij}^a y^j - g^{aj}C_{ijk}\gamma_{lm}^k y^l y^m,$$

where $C_{ijk} = \partial_{y^k} g_{ij}/2$ is the Cartan tensor of $L$, which measures how far $g$ is from being pseudo-Riemannian.

Next our aim is to select an anisotropic Levi-Civita connection for a pseudo-Finsler metric $L$, revisiting the role of the Chern connection.

**Theorem 4.** *Given a pseudo-Finsler metric $L$ and being $g$ its fundamental tensor, there exists a unique anisotropic conection $\nabla$ that is torsion-free and such that $\nabla g = 0$. It is characterized by the Koszul-type formula*

$$\begin{aligned}
2g_v(\nabla_X^V Y, Z) = {} & (X(g_V(Y,Z)) - Z(g_V(X,Y)) + Y(g_V(X,Z)))\,(\pi(v)) \\
& + g_v([X,Y],Z) + g_v([Z,X],Y) - g_v([Y,Z],X) \\
& - 2C_v(Y,Z,\nabla_X^V V) - 2C_v(Z,X,\nabla_Y^V V) + 2C_v(X,Y,\nabla_Z^V V),
\end{aligned} \tag{25}$$

*where $v \in A$, $X,Y,Z \in \mathfrak{X}(M)$, $V \in \mathfrak{X}(\Omega)$ is any local extension of $v$ and $C$ is the Cartan tensor defined above. Moreover, its Christoffel symbols are*

$$\Gamma_{ij}^a = \frac{1}{2} g^{ak} \left( \frac{\delta g_{ki}}{\delta x^j} + \frac{\delta g_{kj}}{\delta x^i} - \frac{\delta g_{ij}}{\delta x^k} \right) \tag{26}$$

*(the $\delta_j$ are the ones of the associated $\nu^\nabla$). Such a unique $\nabla$ is called the* Levi-Civita anisotropic connection *of $L$ and the corresponding vertically trivial linear connection is the* Chern connection.

*Proof.* Taking into account that $\dot{\partial} g = 2C$, it follows that

$$(\nabla_X g)_v(Y,Z) = X(g_V(Y,Z))(\pi(v)) + g_v(\nabla_X^V Y, Z) + g_v(Y, \nabla_X^V Z) - 2C_v(Y,Z,\nabla_X^v V)$$

and using that $\nabla g = 0$, as well as the above formula permuting $X,Y,Z$, one gets (25). To get (26), observe that $\nabla g = 0$ in coordinates means

$$\delta_k g_{ij} - \Gamma_{ki}^l g_{lj} - \Gamma_{kj}^l g_{il} = 0,$$

which is equivalent to the structure equations of the Chern connection (see [2]), and therefore its Chritoffel symbols coincide with those of Chern's as well as its vertically trivial associated connection.  $\qquad\square$

*Remark 8.* We have seen that there are two distinguished anisotropic connections associated with a pseudo-Finsler metric, the Berwald connection (see Remark 7) and the Levi-Civita–Chern connection (see Theorem 4). The difference between them is a tensor $\mathcal{L}^\flat$ metrically equivalent to the Landsberg tensor of $L$, which satisfies $\mathcal{L}_v^\flat(v,u) = 0$ for all $v \in A$ and $u \in T_{\pi(v)}M$ (see [24]). This property is essential as it guarantees that both connections have the same associated non-linear connection $\overset{o}{\nu}$. Indeed, the anisotropic connections differing in a symmetric tensor $Q$ with this property from the Chern or Berwald connections are exactly the torsion-free anisotropic connections having $\overset{o}{\nu}$ as their associated non-linear connection. The properties of this family of connections as a tool for the study of pseudo-Finsler metrics were collected in [11], where they are referred to as the *distinguished connections*.

### 6.2   The Finslerian linear connections

The linear connections associated with a pseudo-Finsler metric $L$ live in the bundle $\mathcal{V}A \to A$. As we have seen, $L$ determines a geodesic spray $G$ and a privileged *metric nonlinear connection* $\overset{o}{\nu}$, which will play the role of the prescribed auxiliary connection seen in §5. Then, Proposition 2 and Theorem 3 are applicable. As a summary, one has:

1. The linear connections $\nabla^*$ used in pseudo-Finsler geometry are defined in the vector bundle $\mathcal{V}A \to A$ and they are homothetically invariant.
2. Such a $\nabla^*$ can be specified by means of the Christoffel symbols with respect to the metric nonlinear connection $\overset{o}{\nu}$ (Proposition 2), namely: $H_{ij}^k$, which are 0-homogeneous in $y$ and $V_{ij}^k$, which are (-1)-homogeneous.
3. The vertically trivial $\nabla^*$'s are in natural correpondence with the homogeneous $A$-anisotropic connections on $M$ (Theorem 3). Using $\overset{o}{\nu}$, the non-vertically trivial $\nabla^*$'s project onto the vertically trivial ones (Remark 6).
4. The most frequent choices of $\nabla^*$ in Finsler Geometry have the following horizontal and vertical parts:
   - Berwald and Hashiguchi: $H_{ij}^k = N_{i;j}^k := \dot{\partial}_j N_i^k$, where the $N_i^j$'s come from $\overset{o}{\nu}$. Berwald is vertically trivial and Hashiguchi has $V_{ij}^k = C_{ij}^k$.
   - Chern-Rund and Cartan: $H_{ij}^k = \mathit{\Gamma}_{ij}^k$ as in (26). Chern-Rund is vertically trivial and Cartan has $V_{ij}^k = C_{ij}^k$.

## 7   Parallel transport and anisotropic connections

Next, let us go back to our observers' viewpoint at §2.2 to introduce the parallel transport and show how the anisotropic connection can be recovered from it.

### 7.1   Observers and parallel transport

The most natural way to compare vectors in different tangent spaces of a manifold is by making use of parallel transport along a curve. Depending on what we want to study, this parallel transport should preserve certain properties of vectors. In general, we cannot ensure the preservation of the indicatrix of a pseudo-Finsler metric by a linear map, because the indicatrices at different points may be not linearly equivalent. The best one that could try to do is to require the preservation of its best approximation by a scalar product, namely, $g_v$. As there is a dependence on $v$, we will need a different parallel transport for every $v \in A$. Summing up, the Christoffel symbols will depend also on the direction, so the covariant derivative along a curve $\gamma : I \to M$ needs a reference vector field $W \in \mathfrak{X}(\gamma)$:

$$D_\gamma^W : \mathfrak{X}(\gamma) \to \mathfrak{X}(\gamma),$$

where $\mathfrak{X}(\gamma)$ will denote the space of smooth vector fields along $\gamma$. Moreover, we will assume that the dependence on $W$ is pointwise, in the sense that at the instant $t_0$, $D_\gamma^W$ depends only on $W(t_0)$. Thinking about what happens in a Finsler spacetime, where all the computations depend on the observer, we will make first the parallel transport of the observer along $\gamma$ by searching for a vector field $V$ such that

$$D_\gamma^V V = 0,$$

and it is natural to require that $L(V) = 1$. Finally, we are ready to parallel transport vectors along $\gamma$ using parallel vector fields $Z$, defined by

$$D_\gamma^V Z = 0.$$

As we will see later, if we require this parallel transport to preserve also the metric $g_v$ of the restspace, then the covariant derivative comes from the Levi-Civita–Chern anisotropic connection. Geodesics are recovered as the curves with autoparallel velocity, namely

$$D_\gamma^{\dot\gamma} \dot\gamma = 0.$$

In particular, $L(\dot\gamma) = const.$

Of course, when we speak about a covariant derivative, we are assuming that it satisfies the natural properties of a derivative. Let us put this on rigorous basis. In the following, given a smooth curve $\gamma : [a, b] \to M$, $\mathcal{F}(I)$ will denote the smooth real functions defined on $I = [a, b]$. Recall that $A$ denotes a conic open subset of $TM$, $\pi : TM \to M$ is the natural projection and $\gamma^*(TM)$ is the pullbacked of this bundle by means of the curve $\gamma : [a, b] \to M$.

**Definition 10.** *An anisotropic covariant derivation $D_\gamma$ in $A$ along a curve $\gamma : [a, b] \to M$ is a map*

$$D_\gamma : \gamma^*(A) \times \mathfrak{X}(\gamma) \to TM, \qquad (v, X) \mapsto D_\gamma^v X \in T_{\pi(v)} M$$

*with a smooth dependence on $v$, such that if $\pi(v) = \gamma(t_0)$ with $t_0 \in [a, b]$,*

*(i)* $D_\gamma^v (X + Y) = D_\gamma^v X + D_\gamma^v Y \ \forall \ X, Y \in \mathfrak{X}(\gamma)$,
*(ii)* $D_\gamma^v (fX) = \frac{df}{dt}(t_0)X(t_0) + f(t_0)D_\gamma^v X \ \forall \ f \in \mathcal{F}(I), \ X \in \mathfrak{X}(\gamma)$.

*Remark 9.* Let $\pi_\gamma : \gamma^*(A) \to [a, b]$ be the pullback fibered manifold obtained from $A \to M$ by $\gamma : [a, b] \to M$. The formal similarity of our definition of $D_\gamma$ with Definition 4 for anisotropic connections can be stressed by redefining: (a) the domain of $D_\gamma$ as $\mathfrak{X}(\gamma)$ and (b) its codomain as the sections of the pullback bundle $\pi_\gamma^*(TM) \to \gamma^*(A)$ obtained from $\gamma^*(TM) \to [a, b]$ by $\pi_\gamma : \gamma^*(A) \to [a, b]$.

The notion of anisotropic connection gathers the information of the covariant derivatives along different curves. In fact, there is a unique covariant derivative along curves determined by an anisotropic connection (see [11, Prop. 2.7]).

**Proposition 4.** *Given a smooth curve $\gamma : [a, b] \to M$, an $A$-anisotropic connection $\nabla$ determines an induced $A$-anisotropic covariant derivative along $\gamma$ with the following property: if $X \in \mathfrak{X}(M)$, then $D_\gamma^v(X_\gamma) = \nabla_{\dot\gamma}^v X$, where $X_\gamma$ is the vector field in $\mathfrak{X}(\gamma)$ defined as $X_\gamma(t) = X_{\gamma(t)} \ \forall t \in [a, b]$.*

Indeed, given local coordinates $(U, x)$, we can express this covariant derivative in terms of the Christoffel symbols of the $A$-anisotropic connection $\nabla$, which are defined as the functions $\Gamma_{ij}^k : T\Omega \cap A \to \mathbb{R}$ determined by

$$\nabla_{\partial_i}^v \partial_j = \Gamma_{ij}^k(v)\, \partial_k|_{\pi(v)}.$$

It is easy to check that if $X = X^k \partial_k$, then

$$D_\gamma^W X = (\dot{X}^i + (\Gamma_{jk}^i \circ W)\dot{\gamma}^j X^k)\partial_i. \tag{27}$$

This provides coordinate expressions for the covariant derivative. Moreover, given a curve $\gamma : [a, b] \to M$, if one fixes the reference vector field $W \in \mathfrak{X}(\gamma)$ and $t_1, t_2 \in [a, b]$, then it is possible to define a parallel transport

$$P_{t_1, t_2}^W : T_{\gamma(t_1)}M \to T_{\gamma(t_2)}M$$

in such a way that $P_{t_1, t_2}^W(z) = Z(t_2)$, being $Z \in \mathfrak{X}(\gamma)$ such that $D_\gamma^W Z = 0$ and $Z(t_1) = z$.

*Remark 10.* The parallel transport $P_{t_1, t_2}^W$ shares all the natural properties of the case of affine connections (it is a well-defined linear isomorphism, naturally invariant under reparametrizations of $\gamma$ including the reversal of the sign). Indeed, as the value of the Christoffel symbols is determined by $W$, which is fixed, (27) yields an equation for the transport of the same type as in the affine case.

There is a different type of parallel transport which is not always well-defined, namely, when the goal is to find a vector field $V$ along $\gamma$ such that $D_\gamma^V V = 0$. The existence of this parallel transport is not guaranteed along the whole curve unless we have some control on the Christoffel symbols.

**Definition 11.** *A smooth curve $\gamma : [a, b] \to M$ is parallel transport complete if for every $v \in A \cap T_{\gamma(a)}M$, there exists a (unique) $A$-admissible vector field $V \in \mathfrak{X}(\gamma)$ such that $D_\gamma^V V = 0$ and $V(a) = v$. Here $A$-admissible means that $V(t) \in A$ for all $t \in [a, b]$.*

*Remark 11.* From standard results of ODE's one has that, for every $v \in A \cap T_{\gamma(a)}M$, there exists some $\epsilon > 0$ such that a parallel $V$ as above is well-defined in $[a, a+\epsilon]$. Moreover, all the curves are parallel transport complete in the following two cases: (1) when $A = TM \setminus \mathbf{0}$ and the anisotropic connection is homogeneous, and (2) in the case of a Finsler spacetime $(M, L)$ with a distinguished connection (for the latter, notice that the anisotropic connection is homogeneous and $L(V)$ is constant for any parallelly transported vector $V$, so that it cannot abandon $A$). From now on, we will restrict ourselves to work with curves where this parallel transport is defined everywhere.

Let us define the parallel transports which have a geometric meaning to compare what happens in different points of the manifold.

**Definition 12 (Instantaneous observer's parallel transport).** *Let* $\nabla$ *be an A-anisotropic connection and* $\gamma : [a, b] \to M$ *a parallel transport complete curve. For each* $t_1, t_2 \in [a, b]$, *the* instantaneous observer's parallel transport *is the map*

$$P_{t_1, t_2} : A \cap T_{\gamma(t_1)} M \to A \cap T_{\gamma(t_2)} M$$

*given by* $P_{t_1, t_2}(v) = V(t_2)$, *for* $V \in \mathfrak{X}(\gamma)$ *satisfying* $V(t_1) = v$ *and* $D_\gamma^V V = 0$.

This parallel transport coincides the one obtained from the nonlinear connection which appears in many classical textbooks devoted to Finsler Geometry as [17, Ch. VII], [1, §2.1.6], [24, page 103], [4, §2.1], [7, Def. 1.4] and [26, §7.6]. In some other textbooks there is an additional notion of parallel transport taking as a reference the velocity of the curve (see [24, Def. 7.3.1],[25, §5.3] and [6, Ch. 4]). Recall that we defined instantaneous observers in the setting of Finsler spacetimes as vectors $v \in A$ of unit length $L(v) = 1$. Of course, the constraint on the length is not relevant for the transport. In the case of general anisotropic connections such a restriction makes no sense but we have maintained the name of observers to stress the relativistic geometric intuitions.

**Definition 13 (Parallel transport with respect to an instantaneous observer).** *Let* $\nabla$ *be A-anisotropic connection and* $\gamma : [a, b] \to M$ *a parallel transport complete curve. For each* $t_1, t_2 \in [a, b]$ *and observer* $v \in T_{\gamma(t_1)} M \cap A$, *the* parallel transport with respect to $v$ *is the map*

$$P_{t_1, t_2}^v : T_{\gamma(t_1)} M \to T_{\gamma(t_2)} M$$

*obtained as* $P_{t_1, t_2}^v(w) = W(t_2)$, *where* $W \in \mathfrak{X}(\gamma)$ *satisfies that* $W(t_1) = w$ *and* $D_\gamma^V W = 0$ *with* $V$ *satisfying* $D_\gamma^V V = 0$ *and* $V(t_1) = v$.

Recall that, as $V$ is fixed by $v$, this parallel transport satisfies the natural properties of the transport explained in Remark 10. See [21] for a general treatment of parallelism.

### 7.2 Recovering the anisotropic connection from the transport

First observe that given a smooth curve $\gamma : [a, b] \to M$, we can define the paralell transport of covectors with respect to the instantaneous observer $v \in T_{\gamma(a)} M \cap A$,

$$P_{a,b}^v : T_{\gamma(a)} M^* \to T_{\gamma(b)} M^*,$$

as

$$P_{a,b}^v(\theta)(w) = \theta(P_{b,a}^v(w))$$

for any $\theta \in T_{\gamma(a)} M^*$ and $w \in T_{\gamma(b)} M$, so that a parallel covector on a parallel vector will be constant. Imposing commutativity with the tensor product, this parallel transport can be extended to arbitrary $(r, s)$-tensors.

Next, let us write explicitly such a transport regarding the tensors as multilinear maps. Consider an $A$-anisotropic tensor $T \in \mathcal{T}_s^r(M_A)$ and a curve $\gamma : [a, b] \to M$. Then we can define the parallel transport $P_{t_1,t_2}(T)_v$ for any $t_1, t_2 \in [a, b]$ as the map

$$P_{t_1,t_2}(T)_v : T_{\pi(v)}M^* \times \overbrace{\cdots}^{r} \times T_{\pi(v)}M^* \times T_{\pi(v)}M \times \overbrace{\cdots}^{s} \times T_{\pi(v)}M \to \mathbb{R}$$

given by

$$P_{t_1,t_2}(T)_v(\theta^1, \ldots, \theta^r, v_1, \ldots, v_s)$$
$$= T_{P_{t_2,t_1}(v)}(P^v_{t_2,t_1}(\theta^1), \ldots, P^v_{t_2,t_1}(\theta^r), P^v_{t_2,t_1}(v_1), \ldots, P^v_{t_2,t_1}(v_s)).$$

In particular, we can define a curve of anisotropic tensors in $T_{\pi(v)}M$:

$$P_t(T) = P_{t,a}(T).$$

Our next goal is to compare this parallel transport with the covariant derivative of any $A$-anisotropic tensor, as given in Theorem 1 and formula (11). Namely, being $P_t(T)$ a curve in a vector space, let us relate its natural derivative with $\nabla_{\dot{\gamma}(a)}T$.

**Proposition 5.** *Given an $A$-anisotropic tensor $T \in \mathcal{T}^r_s(M_A)$, an $A$-anisotropic connection $\nabla$ and a regular curve $\gamma : [a, b] \to M$, it holds that*

$$(\nabla_{\dot{\gamma}(a)}T)_v = \frac{d}{dt}P_t(T)_v|_{t=a},$$

*for any $v \in T_{\gamma(a)}M \cap A$.*

*Proof.* Assume first that $r = 0$. Recall that we can compute $(\nabla_{\dot{\gamma}(a)}T)_v(v_1, \ldots, v_s)$ choosing an $A$-admissible extension $V$ of $v$ and arbitrary extensions $X, X_1, \ldots, X_s$ of $\dot{\gamma}(a), v_1, \ldots, v_s$, respectively. In particular, these extensions can be chosen in such a way that $\nabla^V_X V = \nabla^V_X X_j = 0$ along $\gamma$ for all $j = 1, \ldots, s$. With these choices,

$$(\nabla_{\dot{\gamma}(a)}T)_v(v_1, \ldots, v_s) = X_{\gamma(a)}(T_V(X_1, \ldots, X_s))$$
$$= \frac{d}{dt}(T_{V_{\gamma(t)}}((X_1)_{\gamma(t)}, \ldots, (X_s)_{\gamma(t)}))|_{t=a}. \quad (28)$$

Finally, observe that $V_{\gamma(t)} = P_{a,t}(v)$, since $V(t) := V_{\gamma(t)}$ is a parallel observer along $\gamma$ (recall that $\nabla^V_X V = 0$) and $V(a) = V_{\gamma(a)} = v$. Moreover, $(X_i)_{\gamma(t)} = P^v_{a,t}(v_i)$ since $\nabla^V_X X_i = 0$ and $(X_i)_{\gamma(a)} = v_i$. Replacing these quantities in (28), we get

$$(\nabla_{\dot{\gamma}(a)}T)_v(v_1, \ldots, v_s)|_{t=a}$$
$$= \frac{d}{dt}T_{P_{a,t}(v)}(P^v_{a,t}(v_1), \ldots, P^v_{a,t}(v_s))|_{t=a} = \frac{d}{dt}P_t(T)_v(v_1, \ldots, v_s)|_{t=a},$$

as required.

For the general case, observe that given the covectors $\theta^1, \ldots, \theta^r$, it is possible to choose one-forms $\omega^i$ such that $\nabla_X^V \omega^i = 0$ and $(\omega^i)_{\gamma(a)} = \theta^i$. Then, applying the proposition for $r = 0$, it is possible to show that $(\omega^i)_{\gamma(t)} = P_{a,t}(\theta^i)$. The results follows analogously to the case $r = 0$ by computing the covariant derivative with $V$ as a reference vector and $\omega^1, \ldots, \omega^r, X_1, \ldots, X_s$ as above. □

It is worth pointing out that, as only the parallel transport close to $t = a$ is needed for each chosen $v$, the result can be applied even if the curve is not parallel transport complete (see Remark 11).

### 7.3  Levi-Civita–Chern connection of a Finsler spacetime

Let $(M, L)$ be a pseudo-Finsler manifold with $L : A \to \mathrm{R}$. We aim to find an $A$-anisotropic connection that defines a parallel transport which preserves some metric properties. As we explained in §7, the parallel transport of an instantaneous observer should preserve the $L$-length and the parallel transport with respect to an instantaneous observer should preserve the fundamental tensor $g_v$, namely, for any curve $\gamma : [a, b] \to M$ and $v \in A \cap T_{\gamma(a)}M$ such that the parallel transport of $v$ is well-defined along $\gamma$,

$$g_{P_{a,b}(v)}(P_{a,b}^v(u), P_{a,b}^v(w)) = g_v(u, w) \tag{29}$$

for all $u, w \in T_{\pi(v)}M$. Observe that this implies in particular that $L(v) = L(P_{a,b}(v))$, as

$$L(v) = g_v(v, v) = g_{P_{a,b}(v)}(P_{a,b}^v(v), P_{a,b}^v(v)) = L(P_{a,b}^v(v)) = L(P_{a,b}(v)).$$

**Proposition 6.** *Let $(M, L)$ be a pseudo-Finsler manifold. Then its Levi-Civita–Chern $A$-anisotropic connection is the only torsion-free connection with a parallel transport that preserves the fundamental tensor of $L$.*

*Proof.* Observe that by Proposition 5, the fact that the parallel transport of $\nabla$ preserves the fundamental tensor $g$ as in (29) is equivalent to $\nabla g = 0$, and therefore $\nabla$ is the Levi-Civita–Chern connection of $(M, L)$. □

## Acknowledgments

MAJ was partially supported by the project PGC2018-097046-B-I00 funded by MCIN/ AEI /10.13039/501100011033/ FEDER "Una manera de hacer Europa" and Fundación Séneca project with reference 19901/GERM/15. This work is a result of the activity developed within the framework of the Programme in Support of Excellence Groups of the Región de Murcia, Spain, by Fundación Séneca, Science and Technology Agency of the Región de Murcia.

MS and FFV are partially supported by the project MTM2016-78807-C2-1-P funded by MCIN/ AEI /10.13039/501100011033/ FEDER "Una manera de hacer Europa", by the project A-FQM- 494-UGR18 Programa FEDER-Andalucía 2014-2020, Junta de Andalucía) and by the framework of IMAG-María de Maeztu grant CEX2020-001105-M funded byMCIN/AEI/ 10.13039/50110001103. FFV is partially supported also by an FPU grant (Formación de Profesorado Universitario) from the Spanish Ministry of Universities.